\documentstyle{amsppt}
\NoRunningHeads
%\NoPageNumbers
%\pagewidth{5in}
%\pageheight{7.8in}
%\Monograph
%\input amsppt.sty
\magnification=\magstep1
\hyphenation{co-deter-min-ant co-deter-min-ants pa-ra-met-rised
pre-print pro-pa-gat-ing pro-pa-gate
fel-low-ship Cox-et-er dis-trib-ut-ive}
\def\leaderfill{\leaders\hbox to 1em{\hss.\hss}\hfill}
\def\A{{\Cal A}}

\def\H{{\Cal H}}

\def\idest{i.e.\ }

\def\d{{\delta}}

\def\w{{\omega}}

\def\dt{{\Bbb D \Bbb T}}
\def\ugen{e}
\def\boxit#1{\vbox{\hrule\hbox{\vrule \kern3pt
\vbox{\kern3pt\hbox{#1}\kern3pt}\kern3pt\vrule}\hrule}}
\def\rabbit{\vbox{\hbox{\kern0pt
\vbox{\kern0pt{\hbox{---}}\kern3.5pt}}}}

\def\tableau#1{
        \hbox {
                \hskip -10pt plus0pt minus0pt
                \raise\baselineskip\hbox{
                \offinterlineskip
                \hbox{#1}}
                \hskip0.25em
        }
}

\def\tabCol#1{
\hbox{\vtop{\hrule
\halign{\strut\vrule\hskip0.5em##\hskip0.5em\hfill\vrule\cr\lower0pt
\hbox\bgroup$#1$\egroup \cr}
\hrule
} } \hskip -10.5pt plus0pt minus0pt}

\def\CR{
        $\egroup\cr
        \noalign{\hrule}
        \lower0pt\hbox\bgroup$
}

% Set up the map arrows for commutative diagrams.
%\def\mapright#1{\smash{
%     \mathop{\longrightarrow}\limits^{#1}}}

%Set up macro for commutative diagrams etc. (see Ex. 18.46 in TeXbook)

\topmatter
\title
Generalized Temperley--Lieb algebras and decorated tangles
\endtitle

\author R.M. Green\endauthor
\affil 
Department of Mathematics and Statistics\\
Lancaster University\\
Lancaster LA1 4YF\\
England\\
{\it  E-mail:} r.m.green\@lancs.ac.uk
\vskip 1in
\endaffil

\abstract
We give presentations, by means of diagrammatic generators and
relations, of the analogues of the Temperley--Lieb algebras associated
as Hecke algebra quotients to Coxeter graphs of type $B$ and $D$.
This generalizes Kauffman's diagram calculus for the Temperley--Lieb
algebra.
\endabstract

\endtopmatter

\centerline{\bf To appear in the Journal of Knot Theory and its Ramifications}

\head Introduction \endhead

The Temperley--Lieb algebra is a certain finite dimensional
associative algebra which first arose in \cite{{\bf 14}} in the context
of Potts models in statistical mechanics.  As well as having
applications to physics, the algebra also appears in the framework of 
knot theory, where it is closely related to the Jones polynomial and
isotopy invariants of links.  This relationship is explained in
\cite{{\bf 10}}, where it is shown that the Temperley--Lieb algebra occurs
naturally as a quotient of the Hecke algebra arising from a Coxeter
system of type $A$.

In his thesis, Graham \cite{{\bf 6}} generalized this realization of the
Temperley--Lieb algebra as a Hecke algebra quotient to the case of a
Coxeter system of arbitrary type.  These Hecke algebra quotients are
the eponymous ``generalized Temperley--Lieb algebras''.
Graham classified the finite dimensional generalized Temperley--Lieb
algebras into seven infinite families: $A$, $B$, $D$, $E$,
$F$, $H$ and $I$, where the family of type $A$ gives the original
Temperley--Lieb algebras.

In this paper, we give presentations, by means of diagrammatic
generators and relations, of the generalized Temperley--Lieb algebras
of types $B$ and $D$, building on the work of tom Dieck \cite{{\bf 2}}
and that of Martin and Saleur \cite{{\bf 13}}.  
The motivation behind this is that Kauffman's
pictorial formation of the type $A$ algebras 
\cite{{\bf 11}} has been of great value
when it comes to understanding otherwise purely abstract algebraic
computations such as the representation theory and cellular structures
(in the sense of \cite{{\bf 7}}).  The algebra of type $A$ is also of
great value in knot theory \cite{{\bf 10}}, and is a central object in
the theory of quantum groups.  This is because the type $A$ Artin
group is both the ``standard'' topological braid group, and a
deformation of the symmetric group.  It is possible, although not yet
clear, that the generalized braid groups may be of some importance in
these areas.

\head 1. Generalized Temperley--Lieb algebras of types $B$ and $D$ \endhead

Our main objects of study are the generalized Temperley--Lieb algebras
arising from Coxeter systems of types $B$ and $D$.  These may readily
be described in terms of generators and relations, as we now show.

The information required to define the algebras is encoded in the
relevant Dynkin diagrams.  The Dynkin diagrams of types $B_n$ and
$D_n$ are numbered as in Figures 1 and 2 respectively.

\topcaption{Figure 1} Dynkin diagram of type $B_n$ \endcaption
\centerline{
\hbox to 3.138in{
\vbox to 0.361in{\vfill
        \includegraphics{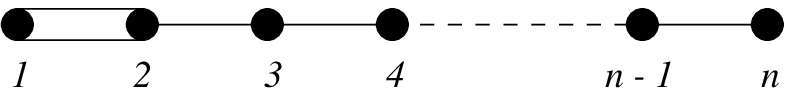}
}
\hfill}
}

\topcaption{Figure 2} Dynkin diagram of type $D_n$ \endcaption
\centerline{
\hbox to 3.388in{
\vbox to 0.916in{\vfill
        \includegraphics{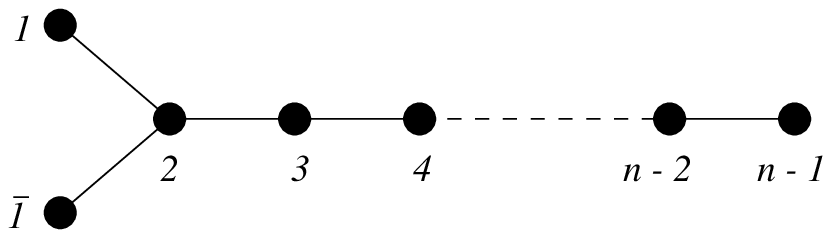}
}
\hfill}
}

We generate the Hecke algebra (as in \cite{{\bf 9}, \S7})
by its Kazhdan--Lusztig basis elements $$
B_s = C'_s := v^{-1} T_1 + v^{-1}T_s
,$$ where $v$ is an indeterminate satisfying $v^2 = q$ and $s$ is a
generating involution of the Coxeter group.

\proclaim{Lemma 1.1}
The Hecke algebra $\H(B_n)$ is generated as an algebra with identity by
the set $\{B_1, B_2, \ldots, B_n \}$ and defining relations $$\eqalignno{
B_s^2 &= (v + v^{-1}) B_s, & (1)\cr
B_s B_t &= B_t B_s \quad \text{ if $|s - t| > 1$}, & (2)\cr
B_s B_t B_s - B_s &= B_t B_s B_t - B_t \quad \text{ if $|s - t| = 1$
and $\{s, t\} \ne \{1, 2\}$}, & (3)\cr
B_s B_t B_s B_t - 2 B_s B_t &= B_t B_s B_t B_s - 2 B_t B_s \quad 
\text{ if $\{s, t\} = \{1, 2\}$.} & (4)\cr
}$$
\endproclaim

\proclaim{Lemma 1.2}
The Hecke algebra $\H(D_n)$ is generated as an algebra with identity by
the set $\{B_{\bar 1},  B_1, B_2, \ldots, B_{n-1} \}$ and defining 
relations $$\eqalignno{
B_s^2 &= (v + v^{-1}) B_s, & (1)\cr
B_s B_t &= B_t B_s \quad \text{ if $s$ and $t$ are not connected in
the Dynkin diagram}, & (2)\cr
B_s B_t B_s - B_s &= B_t B_s B_t - B_t \quad \text{ otherwise}. & (3)\cr
}$$
\endproclaim

\demo{Proofs of Lemmas 1.1 and 1.2}
A routine calculation shows that these definitions are equivalent to
the usual ones involving the $T_s$ elements.
\qed\enddemo

The generalized Temperley--Lieb algebra $TL(X)$ associated to a
Coxeter system with graph $X$ is the quotient of the Hecke algebra $\H(X)$
obtained by factoring out the ideal $I(X)$ generated by the elements $$
\sum_{w \in \langle s_i, s_j \rangle} T_w
$$ as the pairs $(s_i, s_j)$ run over pairs of adjacent nodes in the
Dynkin diagram.

We denote the image in $TL(X)$ of $B_s \in \H(X)$ by $E_s$.

This construction generalizes the construction of the Temperley--Lieb
algebra from the Hecke algebra of type $A$ (see \cite{{\bf 10}}), and is due
to Graham \cite{{\bf 6}}.

It is now not hard to describe the algebras  
$TL(B_n)$ and $TL(D_n)$ explicitly.  

\proclaim{Proposition 1.3}
Let $n \in {\Bbb N} \geq 2$.  We define the associative, unital algebra
$TL(B_n)$ over the ring $\A = {\Bbb Z}[v, v^{-1}]$ 
via generators $E_1, E_2, \ldots E_n$ and 
relations $$\eqalign{
E_i^2 &= [2] E_i, \cr
E_i E_j &= E_j E_i \text{\quad if \ $|i - j| > 1$},\cr
E_i E_j E_i &= E_i \text{\quad if \ $|i - j| = 1$ \ and \ $i, j > 1$},\cr
E_i E_j E_i E_j &= 2 E_i E_j \text{\quad if \ $\{i, j\} = \{1, 2\}$}.\cr
}$$  Here, $[2] := v + v^{-1}$.
\endproclaim

\demo{Proof}
This follows from Lemma 1.1 once it has been observed that the
extra relations imposed on $\H(B_n)$ to make $TL(B_n)$ are precisely
those which set each side of equations (3) and (4) of Lemma 1.1 to zero.
\qed\enddemo

We can apply the same arguments to type $D_n$, as follows.

\proclaim{Proposition 1.4}
Let $n \in {\Bbb N} \geq 4$.  We define the associative, unital algebra
$TL(D_n)$ over the ring $\A = {\Bbb Z}[v, v^{-1}]$ 
via generators $E_1, E_{\bar 1}, E_2, \ldots E_{n-1}$ and 
relations $$\eqalign{
E_i^2 &= [2] E_i, \cr
E_i E_j &= E_j E_i \text{\quad if $i$ and $j$ are not connected
in the graph},\cr
E_i E_j E_i &= E_i \text{\quad if $i$ and $j$ are connected in the
graph.}\cr
}$$
\endproclaim

\demo{Proof}
This follows from Lemma 1.2 along the lines of the proof of
Proposition 1.3.
\qed\enddemo

\head 2. Decorated tangles \endhead

A convenient way to introduce the diagram calculi relevant to this
paper is by means of the category of ``decorated tangles''.  As well as
being an important tool in this paper, the decorated tangles go on to
play a key r\^ole in the sequel to this paper \cite{{\bf 8}}, which
analyses the structure of the generalized Temperley--Lieb algebras of
type $H$.

We introduce a category, the morphisms of which are certain tangles.
Our set-up tends to follow that of
Freyd and Yetter \cite{{\bf 5}}, who showed how Kauffman's tangle-theoretic
approach to the Temperley--Lieb algebra may be defined in terms of
certain categories.

A tangle is a portion of a knot diagram contained in a rectangle.  The
tangle is incident with the boundary of the rectangle only on the
north and south faces, where it intersects transversely.  The
intersections in the north (respectively, south) face are numbered
consecutively starting with node number $1$ at the western (\idest the
leftmost) end.

Two tangles are equal if there exists an isotopy of the plane carrying
one to the other such that the corresponding faces of the rectangle
are preserved setwise.

We call the edges of the rectangular frame ``faces'' to avoid
confusion with the ``edges'' which are the arcs of the tangle.

For our purposes, it is necessary to extend the notion of a
tangle so that each arc of the tangle may be assigned a nonnegative
integer.  (This is similar to the notion of ``coloured'' tangles in
\cite{{\bf 5}}.)  If an arc is assigned the value $r$, we represent this
pictorially by decorating the arc with $r$ blobs.

A decorated tangle is a crossing-free tangle in which each arc is
assigned a nonnegative integer.  Any arc not exposed to the west face
of the rectangular frame must be assigned the integer $0$.
This means that any decorated tangle consists only of loops and edges,
none of which intersect each other.

\example{Example}
Figure 3 shows a typical example of a decorated tangle.  We will
tend to emphasise the intersections of the tangle with the frame
rather than the frame itself, which is why each node (\idest
intersection point with the frame) is denoted by a disc.
\endexample

\topcaption{Figure 3} A decorated tangle \endcaption
\centerline{
\hbox to 3.027in{
\vbox to 0.888in{\vfill
        \includegraphics{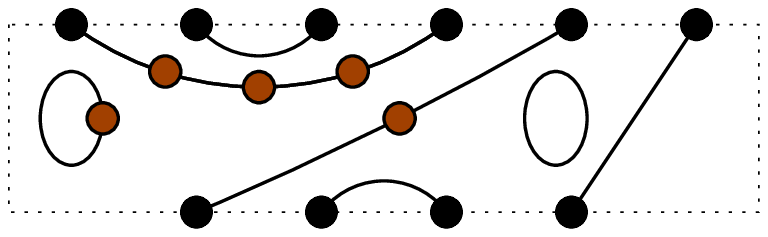}
}
\hfill}
}

The category of decorated tangles, $\dt$, has as its objects the
natural numbers.  The morphisms from $n$ to $m$
are the decorated tangles with $n$ nodes in the north face and $m$ in
the south.  The
source of a morphism is the number of points in the north face of the
bounding rectangle, and the target is the number of points in the
south face.  Composition of morphisms works by concatenation of the
tangles, matching the relevant south and north faces together.

Note that for there to be any morphisms from $n$ to $m$, it is
necessary that $n+m$ be even.  Also notice that the asymmetric
properties of the west face of the rectangle mean that we cannot
introduce the tensor product of two morphisms by the lateral
juxtaposition of diagrams as in \cite{{\bf 5}}.

These category-theoretic definitions allow us to define an algebra of
decorated tangles.
Let $R$ be a commutative ring and let $n$ be a positive integer.  
Then the $R$-algebra $\dt_n$ has as a
free $R$-basis the morphisms from $n$ to $n$, where the multiplication
is given by the composition in $\dt$.

The edges in a tangle $T$ which connect nodes (\idest not the loops)
may be classified
into two kinds: {\it propagating edges}, which link a node in the north
face with a node in the south face, and {\it non-propagating} edges, which
link two nodes in the north face or two nodes in the south face.

It is convenient to define certain named tangles,
$\ugen_{\bar 1}, \ugen_1, \ugen_2, \ldots, \ugen_{n-1}$ and $e$, in
the algebra $\dt_n$.

The tangle $\ugen_i$ is defined as follows.
There is an edge connecting nodes $i$ and
$i+1$ in the north face to each other, and the same for the south
face.  For other nodes $k \ne i, i+1$, node $k$ in the north face is
connected to node $k$ in the south face.  There are no decorated edges
and no loops.

The tangle $\ugen_{\bar 1}$ is obtained from $\ugen_1$ by adding
decorations to the two non-pro\-pa\-gat\-ing edges.

The tangle $e$ has no loops and all its edges are propagating.  There
is one decorated edge, namely the one joining node 1 in the north face
to node 1 in the south face.

\example{Example}
In the case $n = 6$, the tangles $e$, $\ugen_2$ and $\ugen_{\bar 1}$ are as
shown in Figures 4, 5 and 6 respectively.
\endexample

\topcaption{Figure 4} The tangle $e$ in $\dt_6$ \endcaption
\centerline{
\hbox to 2.652in{
\vbox to 0.888in{\vfill
        \includegraphics{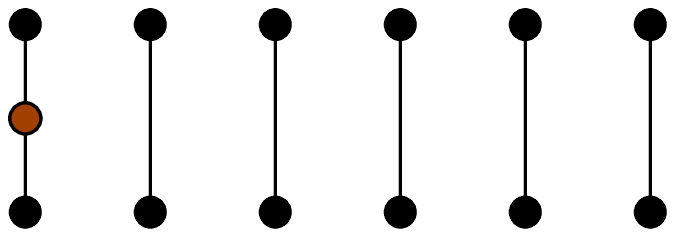}
}
\hfill}
}

\topcaption{Figure 5} The tangle $\ugen_2$ in $\dt_6$ \endcaption
\centerline{
\hbox to 2.638in{
\vbox to 0.888in{\vfill
        \includegraphics{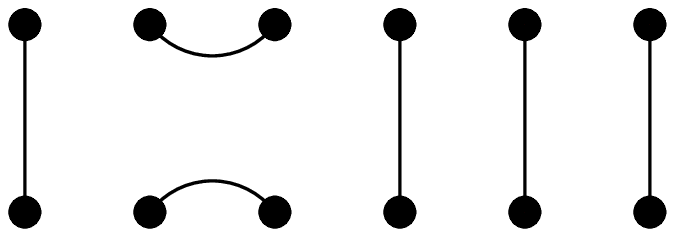}
}
\hfill}
}

\topcaption{Figure 6} The tangle $\ugen_{\bar 1}$ in $\dt_6$ \endcaption
\centerline{
\hbox to 2.638in{
\vbox to 0.888in{\vfill
        \includegraphics{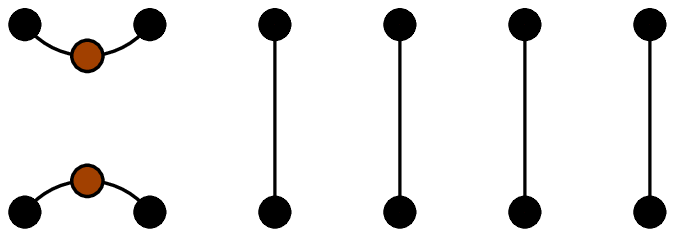}
}
\hfill}
}

Note that if $\ugen_i$ and $\ugen_j$ are such that $i, j \ne {\bar
1}$, then the
relations $\ugen_i \ugen_j = \ugen_j \ugen_i$
(if $|i - j| > 1$) and 
$\ugen_i \ugen_j \ugen_i = \ugen_i$ (if $|i - j| = 1$) hold in $\dt_n$.

\head 3. Review of results in type $A$ \endhead

We now recall Kauffman's tangle-theoretic approach to the
Temperley--Lieb algebra $TL_n$ in terms of the algebra
$\dt_n$ of decorated tangles.  Proofs may be found in \cite{{\bf 12}}.

\proclaim{Theorem 3.1 (Kauffman)}
Let $\d$ be an indeterminate.
Consider the subalgebra of $\dt_n [\d]$ generated by the elements $\ugen_1,
\ugen_2, \ldots, \ugen_{n-1}$ (but not $\ugen_{\bar 1}$).
Let $TL_n$ be the quotient of this subalgebra by the relation

\centerline{
\hbox to 0.916in{
\vbox to 0.277in{\vfill
        \includegraphics{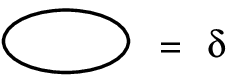}
}
\hfill}
}

The algebra $TL_n$ is the Temperley--Lieb
algebra, which has as a basis the set of all undecorated elements of
$\dt_n$ with no loops.  It is given by generators $\ugen_1,
\ugen_2, \ldots \ugen_{n-1}$ and defining relations $$\eqalign{
\ugen_i ^2 &= \d \ugen_i,\cr
\ugen_i \ugen_{i\pm 1} \ugen_i &= \ugen_i,\cr
\ugen_i \ugen_j &= \ugen_j \ugen_i \quad \text{ if $|i - j| > 1$}.\cr
}$$
\endproclaim

What the relation involving the loop
means is that each occurrence of an undecorated loop
is removed, and the resulting tangle element is
multiplied by the indeterminate $\d$ to compensate.

It should be noted that no decorated edges or loops can arise, since
the generators $e$ and $\ugen_{\bar 1}$ are not involved.

The relation
$\ugen_1^2 = \d \ugen_1$ is illustrated in Figure 7.

\topcaption{Figure 7}
The relation $\ugen_1^2 = \d \ugen_1$ in $TL_3$ \endcaption
\centerline{
\hbox to 3.138in{
\vbox to 1.638in{\vfill
        \includegraphics{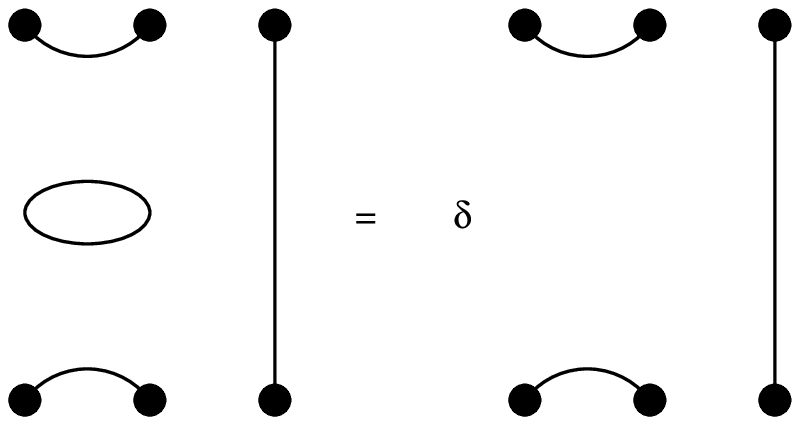}
}
\hfill}
}

The rank of $TL_n$ is well-known to be equal to the Catalan number $$
C(n) := {1 \over {n+1}} {{2n} \choose n}.$$  A basis may be described
in terms of ``reduced words'' in the algebra generators $\ugen_i$.
A reduced word for $TL_n$ is a monomial in the generators
$\{\ugen_1, \ugen_2, \ldots, \ugen_{n-1}\}$ of minimal length.  In
other words, any sequence of relations in Theorem 3.1 which can be applied to
the monomial consists only of applications of commutation
relations.

It is clear that any word in the generators is of the form 
$\d^{a}$ times a reduced word for some integer $a$,
simply by repeated application of the relations.  Thus the
reduced words form a spanning set for the algebra.  In fact, it is
well-known that after discarding repeats, the reduced words give the
same basis of $TL_n$ as the diagrams in Theorem 3.1.

Reduced words for $TL_n$ also have the following important property.

\proclaim{Lemma 3.2}
Let $\ugen_w = \ugen_{i_1} \ugen_{i_2} \cdots \ugen_{i_r}$ be a reduced word in $TL_n$.
Define $a := \min \{ i_1, i_2, \ldots, i_r \}$ and
$b := \max \{ i_1, i_2, \ldots, i_r \}$.  Then there is exactly one
occurrence of $\ugen_a$ in $\ugen_w$, and there is exactly one occurrence of
$\ugen_b$ in $\ugen_w$.
\endproclaim

\demo{Proof}
This is a special case of \cite{{\bf 4}, Lemma 4.3.5}.
\qed\enddemo

Another important property of $TL_n$ is that the number of occurrences
of $\ugen_i$ in a reduced word may be found by inspection of the
corresponding diagram.

\proclaim{Lemma 3.3}
Let $D$ be a basis diagram for the Temperley--Lieb algebra $TL_n$.  Assume the
rectangular frame is drawn with nodes in positions $\{0, 1\} \times
\{1, 2, \ldots, n\}$.  Assume the diagram $D$ is drawn so that the total
number, $2 \ell(D)$, 
of intersections of the associated link with the set of lines
$x = k + 1/2$ (as $k$ runs from $1$ to $n-1$) is minimal.

Let $\ugen_w$ be a reduced monomial in the generators $\{\ugen_1, \ldots,
\ugen_{n-1}\}$ which is equal to $D$.  Then the number of occurrences of
$\ugen_i$ in $\ugen_w$ is half the number of intersections of $D$ with the
line $x = i + 1/2$, and the length of $\ugen_w$ is $\ell(D)$.
\endproclaim

\demo{Proof}
This is a consequence of \cite{{\bf 4}, Lemma 4.3.5}.
\qed\enddemo

\head 4. The main results \endhead

We can now state the two main theorems of this paper, which show how
to realize the generalized Temperley--Lieb algebras of types $B$ and
$D$ in terms of decorated tangles.  Most of the rest of the paper will
be devoted to proving these results.

We first deal with type $B$.  In this case, we assume that the base
ring contains $1/2$, so that in
particular, that we are not in the situation of characteristic $2$.

\topcaption{Figure 8} Relations for type $B$ \endcaption
\centerline{
\hbox to 1.222in{
\vbox to 1.916in{\vfill
        \includegraphics{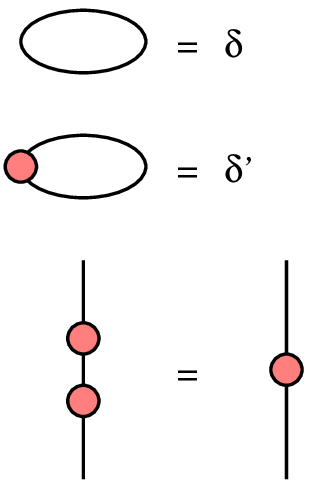}
}
\hfill}
}

\proclaim{Theorem 4.1}
The algebra $TL(B_n)$ arises from $\dt_{n+1}$ as an algebra of 
diagrams via generators $\{\ugen_{\bar 1}, \ugen_2, \ldots, \ugen_n\}$ and
relations shown in Figure 8.

There is a basis for $TL(B_n)$ which is in natural bijection with
elements of $\dt_{n+1}$ which have no loops, at most one decoration on
each edge, and which satisfy one of the following three mutually
exclusive conditions:

\item{\rm ($1$)}
{Node 1 in the north face is joined to node 1 in the south face by an
undecorated edge, and there are no decorated edges.}
\item{\rm ($1'$)}
{Node 1 in the north face is joined to node 1 in the south face by an
decorated edge, but there are no other decorated edges.  Also, there
is at least one non-propagating edge.}
\item{\rm ($2$)}
{The edges emerging from node 1 in the north face and node 1 in the
south face are distinct and both decorated.}

We say that an element of $\dt_{n+1}$ which satisfies these hypotheses 
is $B$-admissible of type $1$, $1'$ or $2$, depending
on which of the three conditions above it satisfies.

This correspondence identifies 
$E_1$ with $2 \ugen_{\bar 1}$ and $E_i$ with $\ugen_i$ for $i > 1$.
\endproclaim

The force of the relations in Figure 8
is firstly to exclude any edge which carries
more than one decoration, and secondly to exclude any loops.  The
third relation in Figure 8
means that all edges and loops may be taken to carry $r$
decorations $(r < 2)$, and the other two relations explain how to remove the
loops.

A simple case by case check verifies that the $B$-admissible diagrams
together with the relations in Figure 8
span an associative algebra (which we will
refer to as $\Cal T(B_n)$).  In particular, the relations are not ambiguous.

The case of type $D$ has a similar overall feel, although there is no
restriction on the characteristic of the base ring.

\topcaption{Figure 9} Relations for type $D$ \endcaption
\centerline{
\hbox to 2.111in{
\vbox to 2.541in{\vfill
        \includegraphics{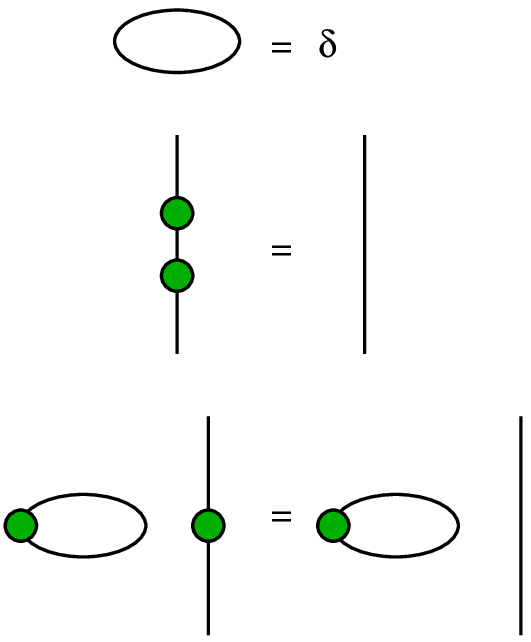}
}
\hfill}
}

\proclaim{Theorem 4.2}
The algebra $TL(D_n)$ arises from $\dt_n$ as an algebra of 
diagrams via generators $\{\ugen_{\bar 1}, \ugen_1, \ugen_2, \ldots, \ugen_{n-1}\}$ and
relations shown in Figure 9.

There is a basis for $TL(D_n)$ which is in natural bijection with
elements of $\dt_n$ which have at most one decoration on each edge or
loop, and which satisfy one of the following two mutually
exclusive conditions:

\item{\rm ($1$)}
{The diagram contains one loop which is decorated, and no other loops
or decorations.  Also, there is at least one non-propagating edge in
the diagram.}
\item{\rm ($2$)}
{The diagram contains no loops and the total number of decorations is
even.}

We say that an element of $\dt_n$ which satisfies these hypotheses 
is $D$-admissible of type $1$ or $2$, depending
on which of the two conditions above it satisfies.

The correspondence identifies
$E_{\bar 1}$ with $\ugen_{\bar 1}$ and $E_i$ with $\ugen_i$ for all
other $i$.
\endproclaim

The third relation in Figure 9 means that any arc loses its
decoration in the presence of a decorated loop.  Using the first and
third relations, all loops may be removed from the image of a diagram
except the last decorated loop, if there is one.  The second relation
ensures that no arc may carry more than one decoration.

Associativity follows by considering triple products of diagrams $D_1
D_2 D_3$, first in the case where the triple product contains a loop
with an odd number of decorations, and then in the other case.

It follows from these observations that
the $D$-admissible diagrams of $\dt_n$ together with the relations in
Figure 9 span an associative algebra (which we will refer to as $\Cal T(D_n)$),
under the diagram multiplication.

\head 5. Combinatorics of decorated tangles \endhead

In order to prove the main results of \S4, we study the combinatorics
of Martin and Saleur's so-called ``blob algebra'', which a
two-parameter version of the algebra studied in \cite{{\bf 2}}, and 
is defined in terms of decorated tangles.  We can associate reduced
words to the
blob algebra, as we did for the Temperley--Lieb algebra in \S3.  The
proofs of theorems 4.1 and 4.2 will be tackled in \S6 by
considering reduced words which satisfy certain additional properties.

The blob algebra $b_n(\d, \d')$ is the algebra which
arises from $\dt_n$ as an algebra of 
diagrams via generators $\{e, \ugen_1, \ugen_2, \ldots, \ugen_{n-1}\}$ and
the relations in Figure 8.  The parameters $\d$ and $\d'$ are indeterminates.

A ``blob diagram'' for $b_n$ is an element of $\dt_n$ which has at most one
decoration on each edge and no loops.

The following result was proved in \cite{{\bf 13}}.

\proclaim{Proposition 5.1 (Martin, Saleur)}
The associative algebra $b_n(\d, \d')$ has as a basis all the blob
diagrams and
has multiplicative structure determined by the relations in Figure 8.
\endproclaim

Notice that the product of two blob diagrams, using
the relations in Figure 8, is a scalar multiple of another one.  The
following lemma is easily verified.

\proclaim{Lemma 5.2}
The following relations hold in $b_n(\d, \d')$: $$\eqalignno{
\ugen_i \ugen_j &= \ugen_j \ugen_i \quad \text{ if $|i - j| > 1$}; & (1)\cr
\ugen_i \ugen_j \ugen_i &= \ugen_i \quad \text{ if $|i - j| = 1$}; & (2)\cr
\ugen_i^2 &= \d \ugen_i; & (3)\cr
e^2 &= e; & (4)\cr
\ugen_1 e \ugen_1 &= \d' \ugen_1; & (5)\cr
\ugen_i e &= e \ugen_i \quad \text{ if $i > 1$}. & (6)\cr
}$$
\endproclaim

A reduced word for $b_n(\d, \d')$ is a monomial in the generators
$\{e, \ugen_1, \ugen_2, \ldots, \ugen_{n-1}\}$ of minimal length.  In
other words,
any sequence of relations in Lemma 5.2 which can be applied to
the monomial consists only of applications of the commutation
relations, (1) and (6).  As in the case of $TL_n$, it is clear that
the reduced words form a spanning set for the algebra.  It also
follows from the results in \cite{{\bf 13}} that the diagram basis for
the blob algebra is the same as the one arising from the reduced words.

By using Lemma 3.2, which concerns extremal generators in reduced
words for $TL_n$, we can show that the occurrences of $e$ and $\ugen_1$ 
in a reduced word for $b_n = b_n(\d, \d')$ alternate.  (Recall
that $e$ commutes with all the generators except $\ugen_1$.)

\proclaim{Lemma 5.3}
Let $\ugen_w$ be a reduced word for $b_n$.  Then there is an occurrence of
$e$ between each pair of occurrences of $\ugen_1$, and an occurrence of
$\ugen_1$ between each pair of occurrences of $e$.
\endproclaim

\demo{Proof}
Suppose there are two occurrences of $e$ in $\ugen_w$.  Then if there is
no occurrence of $\ugen_1$ between them, we can apply relation (6) until
the two occurrences of $e$ are adjacent, and then apply relation (4).
This is a contradiction because $\ugen_w$ is reduced.

Suppose there are two occurrences of $\ugen_1$ in $\ugen_w$, occurring at
positions $c$ and $d$, where $c < d$.  If there is no occurrence of
$e$ between the two occurrences of $\ugen_1$, then $$
\ugen_{i_{c+1}} \ugen_{i_{c+2}} \cdots \ugen_{i_{d-1}}
$$ is a reduced word for $TL_n$ containing two occurrences of $\ugen_1$,
which contradicts Lemma 3.2.  This completes the proof.
\qed\enddemo

We now define two subsets of reduced words for $b_n$: those satisfying
the $B$-condition and those satisfying the $D$-condition.  The reason
for the names is of course that they will be useful in dealing with
the generalized
Temperley--Lieb algebras of types $B$ and $D$, respectively.

Let $\ugen_w$ be a reduced word for $b_n$.

We say $\ugen_w$ satisfies the $B$-condition if one of the following two
conditions holds.
\item{(1)}{Neither $\ugen_1$ nor $e$ occurs in $\ugen_w$.}
\item{(2)}{Both $\ugen_1$ and $e$ occur in $\ugen_w$ but there is no
occurrence of $\ugen_1$ to the left of the leftmost occurrence of $e$ and
there is no occurrence of $\ugen_1$ to the right of the rightmost
occurrence of $e$.}

We say $\ugen_w$ satisfies the $D$-condition if $e$ occurs in $\ugen_w$ an
even number of times (possibly zero).

If a reduced word $\ugen_w$ satisfies the $B$-condition or the
$D$-condition, the following results show that we can essentially
forget about $e$ and restrict our attention to $\ugen_{\bar 1} = e \ugen_1 e$.

\proclaim{Lemma 5.4}
Let $\ugen_w$ be a reduced word satisfying the $B$-condition.  Then $\ugen_w$
is equal to a word $\ugen^{\prime}_w$ in the generators $$
\{\ugen_{\bar 1}, \ugen_2, \ugen_3, \ldots, \ugen_{n-1}\}
.$$
\endproclaim

\demo{Proof}
We describe a procedure for constructing $\ugen^{\prime}_w$.  

First, for each occurrence 
of $e$ in $\ugen_w$ which appears between two occurrences of $\ugen_1$,
replace $e$ by {\sl two} occurrences of $e$.  The $B$-condition now
guarantees that each occurrence of $\ugen_1$ appears between
two occurrences of $e$.

Next, for each occurrence of $\ugen_1$, we can commute the two surrounding
occurrences of $e$ towards the occurrence of 
$\ugen_1$ to form subsequences $e \ugen_1 e$.  This produces a word of the
desired form.
\qed\enddemo

We present some examples from $b_5$
to illustrate the $B$-condition and Lemma 5.4.

\example{Examples}
The words $e \ugen_1$, $\ugen_2 e$, $\ugen_1 \ugen_3$ and $e$ do not satisfy the
$B$-condition, although they are reduced.

The words $\ugen_3 \ugen_2 \ugen_4 \ugen_3$ and $e \ugen_2 \ugen_1 e \ugen_3$ are both reduced and
satisfy the $B$-con\-dition.

The word $\ugen_w = e \ugen_1 e \ugen_2 \ugen_1 e$ satisfies the $B$-condition.  Applying
Lemma 5.4 to $\ugen_w$ first doubles the middle $e$ to form $e \ugen_1 e e
\ugen_2 \ugen_1 e$, and then commutes the generators to form $e \ugen_1 e \ugen_2 e
\ugen_1 e = \ugen_{\bar 1} \ugen_2 \ugen_{\bar 1} = \ugen^{\prime}_w$.
\endexample

The $D$-condition leads to the following property.

\proclaim{Lemma 5.5}
Let $\ugen_w$ be a reduced word satisfying the $D$-condition.  Then $\ugen_w$
is equal to a (not necessarily reduced) word $\ugen^{\prime}_w$ in the
generators $$
\{\ugen_{\bar 1}, \ugen_1, \ugen_2, \ugen_3, \ldots, \ugen_{n-1}\}
.$$
\endproclaim

\demo{Proof}
This is similar to the proof of Lemma 5.4 but slightly simpler.

We use Lemma 5.3 to see that the occurrences of $\ugen_1$ and $e$
alternate in a reduced word $\ugen_w$.  The $D$-condition guarantees that
there is an even number of occurrences of $e$.   Suppose these occur
at positions $i_1, i_2, \ldots, i_{2k}$; we then pair the occurrences
off by twinning the $e$ at position $i_{2j - 1}$ with that at position
$i_{2j}$.  Next we commute each pair of occurrences of $e$ towards the
unique $\ugen_1$ which lies between them to form subexpressions of the
form $e \ugen_1 e$.  This produces a word $\ugen^{\prime}_w$ of the desired form.
\qed\enddemo

We give some examples to illustrate Lemma 5.5.

\example{Examples}
The words $e$ and $\ugen_2 \ugen_1 e \ugen_3$ are reduced but
do not satisfy the $D$-con\-di\-tion.

The words $\ugen_1 \ugen_3$ and $\ugen_w = e \ugen_1 \ugen_2 e \ugen_1$ satisfy the $D$-condition.
Applying Lemma 2.3.4 to $\ugen_w$ pairs off the two occurrences of $e$ and
commutes the rightmost one one place to the left to form $e \ugen_1
e \ugen_2 \ugen_1 = \ugen_{\bar 1} \ugen_2 \ugen_1 = \ugen^{\prime}_w$.
\endexample

The final combinatoric tool needed for the proofs of the main results
is the correspondence between blob diagrams for $b_n$ and diagrams for
$TL_{2n}$ which satisfy a certain symmetry property.
This correspondence produces
tom Dieck's ``symmetric bridges'' \cite{{\bf 2}, \S1}.

Consider a blob diagram $D$ for $b_n$.
Break each decorated edge of $D$ at the decoration, and connect all the loose
endpoints to the west wall in such a way that they do not intersect
each other.  (We will call this the asymmetric representation for $D$.) 
Now consider the diagram union its reflection in the west
wall, which is a diagram for $b_{2n}$ with no decorations.
We will call this the symmetric representation of $D$.

It is not hard to see that this procedure in fact
establishes a bijection between laterally
symmetric diagrams for $TL_{2n}$ and blob diagrams for $b_n$.

\example{Example}
The asymmetric representation of the diagram in Figure 10 is given in
Figure 11, where the west wall is shown explicitly by the dotted line.
The symmetric representation is obtained simply by considering the
west wall as a mirror.
\endexample

\topcaption{Figure 10} A blob diagram for $b_6$ \endcaption
\centerline{
\hbox to 2.638in{
\vbox to 0.888in{\vfill
        \includegraphics{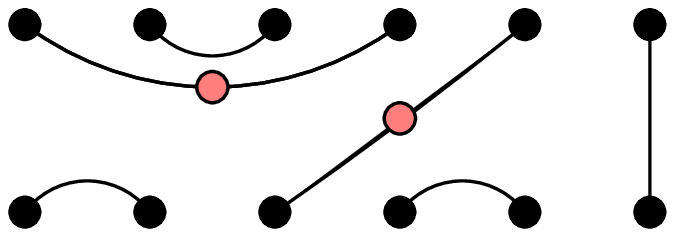}
}
\hfill}
}

\topcaption{Figure 11} Asymmetric representation of Figure 10 \endcaption
\centerline{
\hbox to 2.883in{
\vbox to 0.888in{\vfill
        \includegraphics{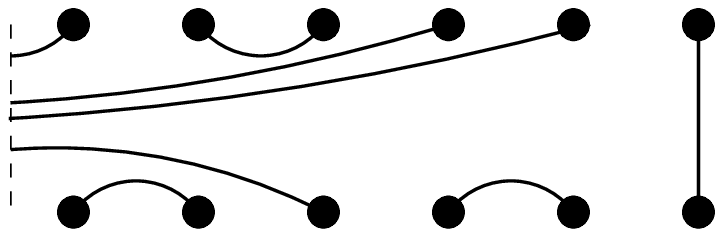}
}
\hfill}
}

From \cite{{\bf 2}, Satz 2.5}, we have

\proclaim{Lemma 5.6}
The left-right symmetric diagrams of $TL_{2n}$ span a subalgebra of
$TL_{2n}$ of dimension ${{2n} \choose n}$
with generators $\{\ugen_1 \ugen_{n-1}, \ugen_2 \ugen_{n - 3}, \ldots,
\ugen_{2n - 1} \ugen_{2n+1}, \ugen_{2n}\}$.
\endproclaim

The next result shows how the algebra of \cite{{\bf 2}} is a special
case of $b_n$.

\proclaim{Lemma 5.7}
If $\d$ is invertible, the subalgebra of $TL_{2n}(\d^2)$ spanned by the 
left-right symmetric diagrams is isomorphic to the algebra
$b_n(\d, 1)$.

The isomorphism may be chosen to identify $\ugen_i$ in $b_n$ with $\ugen_{n+i}
\ugen_{n-i}$ in $TL_{2n}$, and $e$ in $b_n$ with $\ugen_n/\d$ in $TL_{2n}$.

Thus $b_n$ has dimension $${{2n} \choose n}.$$
\endproclaim

\demo{Proof}
The proof of the first two parts 
is simply a matter of checking that the multiplicative action of
the generators is as asserted.

The third part is immediate from the correspondence between the two
bases and Lemma 5.6.
\qed\enddemo

It will be helpful in the proofs of the main theorems to know the 
significance of the total number of decorations in a blob diagram.

\proclaim{Lemma 5.8}
Let $D$ be a blob diagram for $b_n$ and let $\ugen_w$ be a reduced monomial
corresponding to $D$.  Then the number of decorations in $D$ is equal
to the number of occurrences of $e$ in $\ugen_w$.
\endproclaim

\demo{Proof}
Let $D'$ be the symmetric representation of $D$, corresponding to an
element of $TL_{2n}$.  We may assume that $D'$ satisfies the
hypotheses of Lemma 3.3.  Lemma 3.3 shows that the number of
occurrences of $\ugen_n$ in a reduced monomial for $D'$ is equal to the
number of intersections with the line $x = n + 1/2$.  Lemma 5.6
shows that this is equal to the number of occurrences of $e$ in
$\ugen_w$.  Reconstructing $D$ from $D'$, we find that the
number of decorations in $D$ is equal to the number of intersections
of $D'$ with $x = n + 1/2$.  This completes the proof.
\qed\enddemo

\head 6. Proofs of the main results \endhead

We now prove Theorem 4.1.
Until further notice, we replace the parameter $\d$ by $[2]$ and $\d'$
by $[2]/2$.

\proclaim{Lemma 6.1}
There is a homomorphism $\rho_B : TL(B_n) \rightarrow {\Cal T}(B_n)$ 
which takes $E_1$ to $2 \ugen_{\bar 1}$ and $E_i$ to $\ugen_i$ for
$i > 1$.
\endproclaim

\demo{Proof}
This follows by checking that all the relations in Proposition 1.3
hold, which presents no difficulties.
\qed\enddemo

In order to prove that $\rho_B$ is an isomorphism, we need to
enumerate the number of $B$-admissible diagrams of the various types.

\proclaim{Lemma 6.2}
In type $B_n$, the number of $B$-admissible diagrams of type $1$
is $C(n)$, of type $1'$ is $C(n) - 1$ and of type $2$ is ${{2n}
\choose n} - C(n)$.  

This is a total of $(n+2)C(n) - 1$, which is the dimension of $TL(B_n)$.
\endproclaim

\demo{Proof}
The diagrams of type $1$ are in canonical bijection with basis
diagrams for $TL_n$.  The correspondence is given
by removal of the edge joining node 1 in the north face to node 1 in
the south face.  The number of such diagrams is therefore
equal to the Catalan number $$
C(n) := {1 \over {n+1}} {{2n} \choose n}.$$

The case of type $1'$ is similar.  The $-1$ in the formula is due to
the exclusion of the diagram $e$, all of whose edges are propagating (and
one of which is decorated).

Let $D$ be a diagram which is either of type $1'$ or $2$, or
equal to $e$.  Consider the symmetric
representation of $D$ as in \S5; this has the form of a basis
diagram for $TL_{2n+2}$.  Observe that the diagrams which turn up in
this way are precisely the symmetric diagrams in which nodes $n+1$ and
$n+2$ in the north face are joined to each other, and similarly for
the south face.  If we remove these four nodes from the picture, as
well as their associated edges, we have a bijection between the
possibilities for $D$ and the set in the statement of Lemma 5.6.  Thus the
number of diagrams of type $1'$ or $2$ is ${{2n} \choose n} - 1$ as
required: the $-1$ comes from the exclusion of the diagram $e$.

The assertion about the dimension of $TL(B_n)$ follows from \cite{{\bf 3}, \S7.2}.
\qed\enddemo

We now show that the generators given in the statement of Theorem 4.1
do indeed generate ${\Cal T}(B_n)$.

\proclaim{Lemma 6.3}
The algebra ${\Cal T}(B_n)$ is generated by the set 
$\{\ugen_{\bar 1}, \ugen_2, \ldots, \ugen_n\}$.
\endproclaim

\demo{Proof}
Because of Lemma 5.4, we can reduce this problem to showing that any
reduced $B$-admissible diagram is given by a reduced word $\ugen_w$ which
satisfies the $B$-condition.

It is clear that the $B$-admissible diagrams are blob diagrams
and that they span a subalgebra of the blob
algebra $b_n$.  Thus, for any $B$-admissible diagram $D$, there exists
a monomial $\ugen_w$ in the set $\{e, \ugen_1, \ugen_2, \ldots, \ugen_n\}$ which is
equal to a scalar multiple of $D$, since every monomial is a multiple
of a diagram and the monomials form a spanning set.  
By omitting unnecessary terms in $\ugen_w$, we may assume that $\ugen_w$ is
reduced and that the scalar involved is 1 (\idest $\ugen_w = D$).

It remains to show that $\ugen_w$ has the $B$-condition.  If the diagram
$D$ is of type $1$, this follows by Theorem 3.1 because we can
choose $\ugen_w$ in such a way that it avoids all occurrences of $e$ and
$\ugen_1$.  (This uses the embedding of $TL_n$ in $\Cal T(B_n)$ which sends
$\ugen_i \in TL_n$ to $\ugen_{i+1}$.)

If $D$ is of type $1'$ or $2$, then clearly $e D = D e = e D e = D$.
Thus $\ugen_w = e \ugen_w e$.  We now consider the relations in Lemma 5.2
which would need to be applied to the monomial $e \ugen_w e$ in order to
make it reduced.  None of these relations alters the fact that there
can be no occurrence of $\ugen_1$ to the left of the leftmost $e$ or to
the right of the rightmost $e$.  Therefore $\ugen_w$ itself has the
$B$-condition if it is reduced.
\qed\enddemo

The proof of Theorem 4.1 is now complete. \qed

We can now drop the restriction that $\d' = \d/2$.
Instead, we can have $\d = v + v^{-1}$ and $\d' = v' + v^{\prime -1}$, which
corresponds to a quotient of a Hecke algebra $\H(B_n)$ with two
independent parameters $q$ and $Q$, where $q = v^2$ and $Q = v^{\prime
2}$.

We now prove Theorem 4.2.
First, we replace the parameter $\d$ by $[2]$.  This should be regarded
as a change of
notation rather than a restriction, because $v + v^{-1}$ may be
assigned any value if $v$ takes values in an algebraically closed field.

\proclaim{Lemma 6.4}
There is a homomorphism $\rho_D : TL(D_n) \rightarrow {\Cal T}(D_n)$ 
which takes $E_{\bar 1}$ to $\ugen_{\bar 1}$ and $E_i$ to $\ugen_i$ for all
other $i$.
\endproclaim

\demo{Proof}
This follows by checking that all the relations in Proposition 1.4
hold.  The most notable relation is that $\ugen_1 \ugen_{\bar 1} = \ugen_{\bar 1} \ugen_1$.
\qed\enddemo

In order to prove that $\rho_D$ is an isomorphism, we need to
enumerate the number of $D$-admissible diagrams of the various types.

\proclaim{Lemma 6.5}
In type $D_n$, the number of $D$-admissible diagrams of type $1$
is $C(n) - 1$, and the number of type $2$ is ${1 \over 2}{{2n} \choose n}$.

This is a total of $$
\left( {{n + 3} \over 2} \right) C(n) - 1
,$$ which is the dimension of $TL(D_n)$.
\endproclaim

\demo{Proof}
The diagrams of type $1$ are in canonical bijection with the nonidentity basis
diagrams for $TL_n$: the correspondence is given by
removal of the decorated loop.  The number of such diagrams is
therefore $C(n) - 1$.

We argue that the number of diagrams of type $2$ is exactly half the
number of blob diagrams for $b_n$, which we know to be ${{2n} \choose n}$.
Consider the permutation induced on the set of blob diagrams for $b_n$
by the map $\w$ defined as follows.  Let $D$ be such a diagram.  We
define the edge $E$ to be the one connected to the node in
the north-west corner of $D$.
Then $\w(D)$ is obtained from $D$ by toggling the decoration on the
edge $E$, that is, decorating $E$ if $E$ is
undecorated, and removing the decoration from $E$ if $E$ is
decorated.  It is clear that
the orbits of the action of the permutation group generated by $\w$
are all of size 2, and that
exactly one element in each orbit has an even number of decorations.
Thus the number of blob diagrams with an even number of decorations is
exactly half of the total, and the claim follows, completing the proof.

The assertion about the dimension of $TL(D_n)$ follows from \cite{{\bf
3}, \S6.2}.
\qed\enddemo

We now show that the generators given in the statement of Theorem 4.2
do indeed generate ${\Cal T}(D_n)$.

\proclaim{Lemma 6.6}
The algebra ${\Cal T}(D_n)$ is generated by the set $$\{\ugen_{\bar 1},
\ugen_1, \ugen_2, \ldots, \ugen_{n-1}\}.$$
\endproclaim

\demo{Proof}
Because of Lemma 5.5, we can reduce this problem to showing that any
reduced $D$-admissible diagram is given by a reduced word $\ugen_w$ which
satisfies the $D$-condition.

Let $D$ be a $D$-admissible diagram.  If $D$ is of type $1$, then let
$D_A$ be the diagram obtained from $D$ by removing the decorated
loop.  This means $D_A$ is a diagram for $TL_n$, and is equal to
a monomial $\ugen_w$ in the generators $\{\ugen_1, \ugen_2, \ldots, \ugen_{n-1}\}$.
Furthermore, $\ugen_w$ is not trivial since $D$ is not allowed to be the
identity diagram.

Let $a$ be minimal such that $\ugen_a$ occurs in $\ugen_w$.  Define $e'_a$ to
be $$
\ugen_a \ugen_{a-1} \cdots \ugen_2 \ugen_1 \ugen_{\bar 1} \ugen_2 \cdots \ugen_{a-1} \ugen_a
,$$ or $\ugen_1 \ugen_{\bar 1}$ if $a = 1$.  Define $e'_w$ to be the monomial
in the generators for ${\Cal T}(D_n)$ obtained by replacing the
leftmost occurrence of $\ugen_a$ in $\ugen_w$ by $e'_a$.  Then one may easily
check that $e'_w$ gives the diagram $D$.  (Note that we have not
assumed $e'_w$ is reduced.)

Now assume that $D$ is a $D$-admissible diagram of type $2$.  This
means that $D$ has the form of a blob diagram $D_b$ which is equal (by
an argument like that in the proof of Lemma 6.3) to a reduced
monomial $\ugen_w$ in the set $\{e, \ugen_1, \ugen_2, \ldots, \ugen_n\}$.  Since $D$ has
an even number of decorations, Lemma 5.8 shows that $\ugen_w$
contains an even number of occurrences of $e$.  This shows that $\ugen_w$
has the $D$-condition, and can therefore be written as a monomial
$e'_w$ in the generators for ${\Cal T}(D_n)$.

Since $\ugen_w$ is reduced, Lemma 5.8 shows that the number of
decorations in $\ugen_w$ is equal to the number of occurrences of $e$ in
$\ugen_w$.  Therefore the third relation in Figure 8 is never needed in
building up the monomial $\ugen_w$.  The fact that $\ugen_w$ is reduced also
means that the parameters $\d$ and $\d'$ and their associated loops
never appear.  Thus the monomial for ${\Cal T}(D_n)$ which has the
same form as $e'_w$ is equal to $D$ (since none of the diagram
relations in Figure 9 are ever used in multiplying out the
monomial).  The proof now follows.
\qed\enddemo

This completes the proof of Theorem 4.2. \qed

\head 7. Applications \endhead

We conclude by mentioning some of the applications of theorems 4.1 and
4.2.  The details, which are not hard to fill in, are left to the
reader.

A natural idea is to extend the algebra $TL(B_n)$ by adding in the
generator $e$.  This larger algebra can be shown to decompose into a
direct sum of $TL_n$ and $b_n(\d, \d')$, which casts light on the
representation theory of $TL(B_n)$.  

Similarly, by adding the diagram
$G$ (see Figure 12), $TL(D_n)$ may be extended to be isomorphic to a
direct sum of $TL_n$ and an algebra $d_n$ which is half the dimension
of $b_n$.  The algebra $d_n$ can also be constructed from $TL(D_n)$ by
treating any loop carrying an odd number of decorations as zero.

\topcaption{Figure 12} The diagram $G$ for $n = 6$ \endcaption
\centerline{
\hbox to 2.958in{
\vbox to 0.888in{\vfill
        \includegraphics{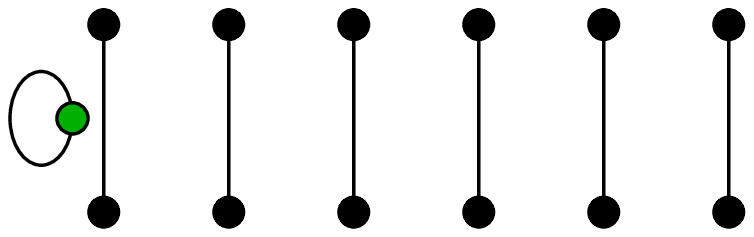}
}
\hfill}
}

Using the results of \cite{{\bf 13}, \S4}, vertex representations of
these algebras may be constructed.  These are representations on a
tensor space $V^{\otimes n}$, where $V$ is of dimension $2$; such a
construction is familiar from statistical mechanics \cite{{\bf 1}}.
The algebra $d_n$ has a particularly natural realisation viewed in
this way as a subalgebra of $b_n$.

One of the important applications of the diagram calculi is the description of
the cellular structures of the algebras $TL(B_n)$ and $TL(D_n)$.  This
uses techniques similar to those in type $A$: a two-sided cell (in
the sense of \cite{{\bf 7}}) consists of a set of diagrams with
particular combinatoric properties.  These properties include the
number and type of propagating edges in the diagram.  Viewed in this
way, the diagrams can be split into a top part and a bottom part,
to form parenthesis diagrams generalizing
those for the Temperley--Lieb algebra \cite{{\bf 15}, \S2}.  The
parenthesis diagrams can then be used to formulate branching rules for
the generically irreducible modules, thus showing for example how such
modules for $TL(B_n)$ restrict to modules for $TL(B_{n-1})$.

It is hoped that there will be further applications of these results
to operator algebras and subfactors.

\head Acknowledgements \endhead
The work for this paper was done while the 
author was supported in part by an E.P.S.R.C. postdoctoral
research assistantship.  The author is grateful to the referee for pointing
out some errors in an earlier version of this paper.

%\leftheadtext{References}
%\rightheadtext{References}

\vfill\eject
\end